\documentclass[11pt,reqno,a4paper]{article}

\usepackage[blocks]{authblk}
\usepackage{amsmath,amsthm,amstext,amscd,amssymb,euscript,url,mathrsfs,euscript,dsfont,calrsfs,mathtools,comment,relsize,xcolor}
\usepackage{hyperref}
\definecolor{unbleu}{rgb}{0.03, 0.15, 0.4}
\hypersetup{
pdfborder = {0 0 0},
colorlinks,
linkcolor=unbleu,
citecolor=unbleu,
urlcolor=unbleu
}

\newcommand{\R}{\mathds R}
\newcommand{\N}{\mathds N}

\newcommand{\E}{\mathds E}
\newcommand{\Zd}{\mathds Z^d}

\newcommand\var{\textup{Var}}
\newcommand{\La}{\ensuremath{\Lambda}}
\newcommand{\la}{\ensuremath{\Lambda}}
\newcommand{\si}{\ensuremath{\sigma}}

\newcommand{\pee}{\ensuremath{\mathbb{P}}}
\newcommand{\loc}{\mathcal{L}}
\def\1{{\mathds 1}}

\newcommand{\dd}{\mathop{}\!\mathrm{d}}
\DeclareMathOperator{\e}{\mathlarger{\mathrm{e}}}

\newtheorem{theorem}{{\small T}{\scriptsize HEOREM}}[section]
\newtheorem{corollary}{{\bf{\small C}{\scriptsize OROLLARY}}}[section]
\newtheorem{proposition}{{\bf{\small P}{\scriptsize ROPOSITION}}}[section]
\newtheorem{lemma}{{\bf{\small L}{\scriptsize EMMA}}}[section]
\newtheorem{remark}{{\bf{\small R}{\scriptsize EMARK}}}[section]
\newtheorem{definition}{{\bf{\small D}{\scriptsize EFINITION}}}[section]

\renewenvironment{proof}[1]
{\noindent{{\bf{\small{ P}{\scriptsize ROOF}}}.}\hspace{0.1cm} #1} {$\;\qed$\newline}

\newcommand{\be}{\begin{equation}}
\newcommand{\ee}{\end{equation}}

\newcommand{\caA}{{\mathcal A}}
\newcommand{\caB}{{\mathcal B}}
\newcommand{\caC}{{\mathscr C}}
\newcommand{\caD}{{\EuScript D}}

\newcommand{\euH}{{\EuScript F}}
\newcommand{\euJ}{{\EuScript J}}
\newcommand{\caF}{{\mathcal F}}
\newcommand{\caG}{{\mathcal G}}

\newcommand{\caL}{{\mathcal L}}

\newcommand{\caS}{{\mathcal S}}

\newcommand{\sipr}{\{\si(t), t\geq 0\}}
\newcommand{\gcb}[1]{\mathrm{GCB}\!\left(#1\right)}
\newcommand{\uvi}[1]{\mathrm{UVB}\!\left(#1\right)}

\begin{document}

\title{Evolution of concentration under\\ lattice spin-flip dynamics}

\author[1]{Jean-Ren\'e Chazottes
\thanks{Email: \texttt{jeanrene@cpht.polytechnique.fr}}
}
\author[1]{Pierre Collet
\thanks{Email: \texttt{pierre.collet@cpht.polytechnique.fr}}
}
\author[2]{Frank Redig
\thanks{Email: \texttt{F.H.J.Redig@tudelft.nl}}
}

\affil[1]{Centre de Physique Th\'eorique, CNRS, Ecole polytechnique, IP Paris, Palaiseau, France}
\affil[2]{Institute of Applied Mathematics, Delft University of Technology, Delft,
%van Mourik Broekmanweg 6, 2628 XE, Delft,
The Netherlands}

\date{Dated: \today}

\maketitle

\begin{abstract}
We consider spin-flip dynamics of Ising lattice spin systems  and study the time evolution of concentration inequalities.
For ``weakly interacting'' dynamics we show that the Gaussian concentration bound is conserved in the course of time and it is satisfied by the unique stationary Gibbs measure.
Next we show that, for a general class of translation-invariant spin-flip dynamics, it is impossible to evolve in finite time from a low-temperature Gibbs state
towards a measure satisfying the Gaussian concentration bound.
Finally, we consider the time evolution of the weaker uniform variance bound, and show that this bound is conserved under a general class of spin-flip dynamics.

\smallskip

\end{abstract}

\newpage

\tableofcontents

%%% SECTION
\section{Introduction}

Concentration inequalities are important tools to understand the fluctuation properties of general observables $f(\si_1, \ldots, \si_n)$ which are functions
of $n$ random variables $(\si_1, \ldots, \si_n)$, where $n$ is large but finite. For bounded random variables which are independent (or weakly dependent)
typically one can obtain so-called Gaussian concentration bounds for the fluctuations of $f(\si_1, \ldots, \si_n)$ about its expectation.
In the context of lattice spin systems, one has, e.g., $\si_i\in\{-1,+1\}$, with $i\in\left[-n,n\right]^d\cap \Zd$, and these random variables are
distributed according to a Gibbs measure. The ``weak dependence'' between them means for instance that we are
in the Dobrushin uniqueness regime, which is for instance the case at ``high enough'' temperature for every finite-range potential, or for low temperature with a ``high enough'' external magnetic
field. In this case a Gaussian concentration bound holds \cite{kuelske}. In contrast, regimes of non-uniqueness are known in which only weaker concentration bounds, such as moment bounds, hold
\cite{cckr}.
In \cite{moles} it is shown that the Gaussian concentration bound implies uniqueness of equilibrium states (translation-invariant Gibbs measures).
In \cite{ccr}, many applications of these concentration bounds are given (speed of
convergence of the empirical measure in the sense  of Kantorovich distance, fluctuation bounds in the Shannon-McMillan-Breiman theorem, fluctuation bounds for the first occurrence of a
pattern, etc).

The Gaussian concentration bound implies volume large deviations for ergodic averages of local observables, i.e., when it holds, the
probability that empirical averages of local observables deviate from their expectation is exponentially small in the volume over which the empirical average is taken. This excludes sub-volume large deviations, which in the context of equilibrium systems implies that the
Gaussian concentration bound cannot hold in a phase transition regime.

%\color{red}
In this paper we are interested in the time evolution of the Gaussian concentration bound under a stochastic evolution.
More precisely we study the following questions in the context of spin-flip dynamics of lattice spin systems:
\begin{enumerate}
\item When started from a probability measure satisfying the Gaussian concentration bound, do we have this bound at later times?
\item When started from a probability measure which does not satisfy the Gaussian concentration bound, can this bound be obtained at finite times?
\end{enumerate}
At the end of the paper we study the same questions for a weaker concentration bound, namely the uniform variance bound.

The study of time-dependent concentration properties of a measure under a stochastic evolution has several motivations.
First, it reveals properties of transient non-equilibrium states, i.e., when one heats up or cools down a system, then what are the concentration properties of the transient states in the course of this process? As mentioned before, the Gaussian concentration bound is a signature
of ``high-temperature'', ``strong uniqueness'' or ``strong mixing''. When cooling or heating a high-temperature system, one can ask whether this signature of
high-temperature behavior is conserved in the course of time, even if one cannot make sense of intermediate temperatures in the course of the evolution, due to possible Gibbs-non-Gibbs transitions \cite{efhr}. Conversely, if one heats a system initially at low temperature,
can the Gaussian concentration bound hold at finite times, i.e., can one obtain this signature of high-temperature behavior in finite time?

Second, semigroups corresponding to stochastic evolution are useful interpolation tools, which give access to properties of measures
which are not available in explicit (e.g. Gibbsian) form. The study of time evolution of concentration properties gives insight in the concentration properties of such measures. An example is e.g. a spin-flip dynamics associated to two different temperatures, where the stationary distribution is an example of a non-equilibrium steady state about which little explicit information is available, as it will generically not be a Gibbs measure (equilibrium state). If one can show the conservation of the Gaussian concentration bound in the course of such a non-equilibrium time evolution, with constants uniformly bounded in time, then one obtains also the Gaussian concentration bound for the non-equilibrium steady state.

Third, the study of time-dependent concentration properties is related to the study of Gibbs-non-Gibbs transitions \cite{efhr}. Here in the regime
where the time-evolved measure is not Gibbs measure, one still would like to obtain some properties of these non-Gibbsian states.
E.g. if one starts a high temperature dynamics from the  low-temperature Ising model with a weak magnetic field, it is known that
one can have Gibbs-non-Gibbs transitions. On the other hand, due to the magnetic field, the initial state satisfies the Gaussian concentration inequality, and therefore if this inequality is conserved in the course of time, one obtains that even in the non-Gibbsian regime, the
measures in the course of time still  satisfy  the Gaussian concentration inequality.
One can also start from the low-temperature Ising model in the phase transition regime and run a high-temperature dynamics.
Then it is also known that in the course of time the Gibbs property is lost, even if the dynamics eventually converges to a high-temperature Gibbs state. It is then interesting so see whether the non-Gibbsian states reached in the course of the time evolution can already have at finite times signatures of the high-temperature behavior of the stationary state, such as the Gaussian concentration bound.
%\color{black}
In the context of time evolution of Gibbs measures, one has generically two scenarios. In the high-temperature regime, i.e.,
high-temperature initial Gibbs measure, and high-temperature dynamics, the time-evolved measure is generically high-temperature Gibbs, and
results of this type are proved via some form of high-temperature (cluster, polymer) expansion, see \cite{efhr}, \cite{mn}. In the regime where the dynamics is high-temperature and the initial measure is
low-temperature, one typically has  Gibbs-non-Gibbs transitions, i.e., after a finite time
the time-evolved measure is no longer a Gibbs measure, and sometimes (e.g. for independent spin-flip dynamics starting
from a low-temperature Ising state with positive small magnetic field) the measure can become Gibbs again.

In the context of time-evolution of concentration inequalities, in \cite{collet3} results so far are restricted to dynamics of diffusive type, in a finite-dimensional context. Here we are interested in the
setting of translation-invariant spin-flip dynamics in infinite
volume, which is precisely the context of Gibbs-non-Gibbs transitions in \cite{efhr}.
Guided by the intuition coming from this context, one expects that a high-temperature dynamics should conserve the Gaussian concentration bound.

We prove this result in the present paper, using the expansion in \cite{mn}, i.e., under the condition that the flip rates are sufficiently close to the rates of an independent spin-flip dynamics.

Next we show that whenever one starts from a low-temperature initial state, i.e., in the non-uniqueness regime, then for any
finite-range spin-flip dynamics, at any later time the distribution cannot satisfy the Gaussian concentration bound. This can be thought of as a result showing that in finite time one cannot obtain
``high-temperature properties'' when initially started from a ``low-temperature state''. This result is shown via an analyticity argument, which
shows that two different initial measures can never coincide in finite time, together with the fact that if a measure satisfies the Gaussian concentration bound, then its
lower relative entropy density with respect to any other translation-invariant measure is strictly positive.
I.e., the existence of two time-evolved measures with zero relative entropy density excludes the possibility that one of them satisfies the Gaussian concentration bound.

%\color{red}
Finally, we show that a weaker concentration bound, the uniform variance inequality, is generically conserved in the course of quasilocal
spin-flip dynamics. This weaker bound which is also valid for pure phases at low temperatures (such as the low-temperature Ising model) implies
that the variance of empirical averages of local observables decays like the inverse of the volume over which the empirical average is taken.
In particular, this excludes divergence of susceptibility, i.e., critical behavior.
Our result implies that in the course of a time evolution started from a non-critical state, no critical state can be obtained.
E.g. if one heats up a low-temperature Ising model, in the limit one obtains a high-temperature state, and in the course of the evolution one never reaches a state which looks like the Ising model at the critical temperature.

The rest of our paper is organized as follows.
In section 2 we introduce some basic context and background on
on Gibbs measures and spin-flip dynamics.
In section 3 we show conservation of the Gaussian concentration bound under a strong high-temperature (or weak interaction) condition.
In section 4 we prove that the Gaussian concentration bound cannot be obtained in finite time if one starts from an initial Gibbs measure in a non-uniqueness (``low-temperature'') regime. In this
section we also prove a non-degeneracy result, based on analyticity, which is of independent interest.
In section 5 we show conservation of the uniform variance inequality for general quasilocal spin-flip dynamics.
%\color{black}

%%% SECTION
\section{Setting: lattice spin systems, Gibbs measures, Markovian dynamics}\label{sec:setting}
%\color{red}
In this section we introduce some basic notation, definition of the Gaussian concentration bounds,  basic concepts about Gibbs measures, spin-flip dynamics and relative entropy.
The expert reader can skip this section, or go over it very quickly.
%\color{black}
We consider the state space of Ising spins on the lattice $\Zd$, i.e., $\Omega=\{-1,1\}^{\Zd}$.
For elements $\si\in\Omega$, called ``spin-configurations'', we denote $\si_i\in \{-1,1\}$ the value of the spin at lattice site $i\in\Zd$.
When we say ``a probability measure $\mu$ on $\Omega$'', we mean a probability measure on the Borel-$\si$-field of $\Omega$, equipped with the standard product of discrete topologies, which
makes $\Omega$ into a compact metric space.
For $\eta\in\Omega$ we denote $\tau_i \eta$ the shifted or translated configuration, defined via $(\tau_i\eta)_j= \eta_{i+j}$.
A function $f:\Omega\to\R$ is called local if it depends only on a finite number of coordinates.
By the Stone-Weierstrass theorem, the set of local functions is dense in the Banach space of continuous functions $\caC(\Omega)$, equipped with the supremum norm.
For $f:\Omega\to\R$ we denote $\tau_i f$ the function defined via $\tau_i f(\eta)= f(\tau_i \eta)$.

For a function $f:\Omega\to\R$ we denote
the discrete gradient
\[
\nabla_{\!i} f (\si)= f(\si^i)- f(\si)
\]
where $\si^i$ denotes the configuration obtained from $\si$ by flipping the symbol at lattice site $i\in\Zd$.
We further denote
\[
\delta_i f=\sup_{\si\in\Omega} \nabla_{\!i} f(\si).
\]
We think of $\delta_i f$ as ``the Lipschitz constant in the coordinate $\si_i$''.
The symbol $\delta f$ means the collection of $\delta_i f , i\in\Zd$, i.e., the ``vector'' of Lipschitz constants.
For $p\geq  1$ we define
\[
\| \delta f\|_p= \left(\sum_{i\in\Zd} (\delta_i f)^p\right)^{\frac{1}{p}}.
\]
%The case $p=2$ will the most important one for our purposes, and therefore
%we abbreviate $\|\delta f\|_2=:\|\delta f\|$.
For a continuous function $f:\Omega\to\R$ and a probability measure $\mu$ on $\Omega$, we will write either $\E_\mu(f)$ or $\int f \dd\mu$ for the integral of $f$ with respect to $\mu$.
We can now define what we mean by a Gaussian concentration bound for a given probability measure on $\Omega$.
\begin{definition}[Gaussian Concentration Bound]
\leavevmode\\
A probability measure $\mu$ on $\Omega$ is said to satisfy the Gaussian concentration bound with constant
$C>0$, abbreviated $\gcb{C}$, if for all continuous $f:\Omega\to\R$ we have
\be\label{gcb}
\E_\mu\left( \e^{f-\E_\mu(f)}\right)\leq \e^{C \|\delta f\|_2^2}.
\ee
\end{definition}
Observe that $\|\delta f\|_2$ is always finite for local functions. Note that a function $f:\Omega\to\R$ is local if and only if there exists a finite subset of $\Zd$ (depending of course on $f$) such that
$\delta_i f=0$ for all $i$ outside of that subset. For non-local continuous functions, inequality \eqref{gcb} is meaningful only when $\|\delta f\|_2<+\infty$.
By a standard argument (exponential Chebyshev inequality applied to $\lambda f$, $\lambda>0$, and then optimization over $\lambda$), the bound \eqref{gcb} implies the ``sub-gaussian'' concentration inequality
\[
\mu\big(f-\E_\mu(f)\geq u\big)\leq \e^{-\frac{u^2}{4C \|\delta f\|_2^2}}
\]
for all $u>0$.
%As mentioned in the introduction, a Gaussian concentration bound has many applications \cite{ccr}.

%%% subSECTION
%\color{red}
\begin{remark}
The Gaussian concentration bound implies in particular ``volume'' large-deviation upper bounds for empirical averages. More precisely, for
a translation-invariant measure satisfying \eqref{gcb}, for a local function $f$, we have
\[
\mu\left(\sum_{x\in \la}\tau_x f-\E_\mu(f)\geq u\right)\leq e^{-|\la| C_f u^2},
\]
with $C_f>0$.
Therefore in the context of Gibbs measures (equilibrium states), it is impossible to have the Gaussian concentration bound
in the non-uniqueness regime.
In this sense, the Gaussian concentration bound can be seen as a signature of ``high-temperature'' or
``weak interaction'' regime. The Gaussian concentration bound is (strictly) weaker than the log-Sobolev inequality, which is the context
of Gibbs measures is known to be equivalent with strong uniqueness conditions \cite{sz}.
\end{remark}
%\color{black}
\subsection{Gibbs measures}

In the context of Gibbs measures, the Gaussian concentration bound is satisfied in the so-called high-temperature regime, and more generally in regimes where the unique Gibbs measure is sufficiently close to a product measure such as the Dobrushin uniqueness regime.
In this subsection we provide some basic background material on Gibbs measures which we need in the sequel. We refer to \cite{geo} for more details and further background.
Let $\caS$ denote the set of finite subsets of $\Zd$. For $\La\subset\Zd$, we denote by $\caF_\la$ the $\si$-field generated by $\{\si_i, i\in \La\}$.
\begin{definition}
%\label{pot}
A uniformly absolutely summable potential  is a map
$U: \caS\times \Omega\to\R$ with the following properties:
\begin{enumerate}
\item
$U(A, \cdot)$ only depends on $\si_i, i\in A$.
\item Uniform absolute summability:
\[
%\label{absum}
\sup_{i\in\Zd} \sum_{\substack{A\in \caS \\ A\ni i}} \,\sup_{\si\in\,\Omega} |U(A, \si)|<+\infty.
\]
\end{enumerate}
\end{definition}
A potential is called translation invariant if $U(A+i, \si)= U(A, \tau_i\si)$ for all $A\in \caS, \si\in \Omega$, $i\in\Zd$.

Given a uniformly absolutely summable potential $U$, and $\La\in \caS$, we denote the finite-volume  Hamiltonian with boundary condition $\eta\in\Omega$:
\[
%\label{finham}
H^\eta_\la(\si_\La)=\sum_{A\,\cap\,\la\not=\emptyset} U(A, \si_\la\eta_{\La^{\!c}})
\]
and the corresponding finite-volume Gibbs measure with boundary condition $\eta$
\[
%\label{finvolgibbs}
\mu^\eta_\la (\si_\La)= \frac{\e^{-H^\eta_\la(\si_\La)}}{Z^\eta_\La}
\]
where $Z^\eta_\La=\sum_{\si_\La\in\Omega_\La} \e^{-H^\eta_\la(\si_\La)}$, the partition function with boundary condition $\eta$, is the normalizing constant (and where $
\Omega_\La$ is the restriction of $\Omega$ to $\La$).
\begin{definition}
Let $U$ be a uniformly absolutely summable potential.
A measure $\mu$ is called a Gibbs measure with potential $U$ if its conditional probabilities satisfy
\[
%\label{dlr}
\mu\big(\si_\la|\caF_{\la^c}\big) (\eta)= \mu^\eta_{\la}(\si_\la)
\]
for all $\la\in\caS$, for all $\sigma$, and for $\mu$-almost every $\eta$.
We will write $\mu\in \caG(U)$ to mean that $\mu$ is a Gibbs measure for $U$.
\end{definition}

We say that $U$ satisfies the strong uniqueness condition if
\be\label{dobu}
c(U):=\sup_{i\in\Zd}\frac12 \sum_{A\ni i} \, (|A|-1)\sup_{\si, \eta\, \in\,\Omega} |\,U(A, \si)- U(A, \eta)|< 1.
\ee
If $U$ satisfies \eqref{dobu} then the set of Gibbs measures $\caG(U)$ is a  singleton (unique Gibbs measure, no phase transition). The condition
\eqref{dobu} implies the well-known Dobrushin uniqueness condition (cf.\ \cite{geo} chapter 8).

If $U$ is translation invariant then $\caG(U)$
contains at least one translation-invariant Gibbs measure.

The following result is a particular case of the main theorem in \cite{kuelske} which states that, under the Dobrushin uniqueness condition, one has the Gaussian concentration bound \eqref{gcb}.
\begin{theorem}[\cite{kuelske}]
\label{kuelthm}
If $U$ satisfies \eqref{dobu} then $\mu\in \caG(U)$ satisfies $\gcb{C}$ with
$C=\frac{1}{2(1-c(U))^2}$.
\end{theorem}
From the proof, one easily infers that also all the finite-volume Gibbs measures $\mu_\la^\eta$ satisfy $\gcb{C}$ whenever
$U$ satisfies \eqref{dobu}, with a constant $C$ that neither depends on the boundary condition $\eta$ nor on the volume $\la$.

%%% subSECTION
\subsection{Relative entropy density and large deviations}

Translation-invariant Gibbs measures with a translation-invariant uniformly absolutely summable potential satisfy a level-3 large deviation principle with the relative entropy density as rate function
\cite[Chapter 15]{geo}.
Let $U$ be a translation-invariant uniformly absolutely summable potential, and $\mu\in \caG(U)$ be a translation-invariant Gibbs measure.
Let $\nu$ be a translation-invariant probability measure on $\Omega$.
The relative entropy density is defined to be the limit
\be\label{relentdef}
h(\nu|\mu)= \lim_{n\to\infty}\frac{\,h_{\la_n}(\nu|\mu)}{|\la_n|}
\ee
with $\la_n= [-n,n]^d\cap \Zd$, $|\la_n|=(2n+1)^d$, and
\[
h_{\la_n}(\nu|\mu)= \sum_{\si_{\la_n}\in\,\Omega_{\La_n}} \nu(\si_{\la_n})\log \frac{\nu(\si_{\la_n})}{\mu(\si_{\la_n})}.
\]
The relative entropy density exists for any $\mu\in \caG(U)$ translation-invariant Gibbs measure, and $\nu$ any translation-invariant probability measure.
Moreover, the relative entropy density is the rate function of the so-called level 3 large deviation principle, i.e.,
in the sense of the large deviation principle, it holds that
\[
%\label{ldp}
\mu\left(\frac1{|\la_n|}\sum_{i\in \la_n}\delta_{\tau_i \si}\approx \nu \right)\asymp \e^{-|\la_n| \, h(\nu|\mu)}.
\]
(This is of course an informal statement where ``$\approx \nu$'' means a neighborhood of $\nu$ in weak topology, and ``$\asymp$'' means asymptotic equivalence after taking
the logarithm and dividing out by $|\la_n|$.)
In general, i.e., if $\mu$ is not a Gibbs measure, the limit defining \eqref{relentdef} might not exist, in that case
we define the lower relative entropy density as
\[
%\label{lowrelent}
h_*(\nu|\mu)=\liminf_{n\to\infty}\frac{\, h_{\la_n}(\nu|\mu)}{|\la_n|} .
\]

The following elementary lemma, which we formulate in the context of a finite set, with a Markov transition matrix,
shows that the relative entropy is decreasing under the action of a Markov kernel.
\begin{lemma}\label{relentlem}
Let $P(x,y)$ be a Markov transition function on a finite set $S$, $x,y\in S$, i.e.,
$P(x,y)\geq 0, \sum_{y\in S} P(x,y)=1$ for all $x\in S$.
Let $\mu,\nu$ be two probability measures on $S$ and let
\[
H(\mu|\nu)= \sum_{x\in S} \mu(x) \log\frac{ \mu(x)}{\nu(x)}
\]
denote their relative entropy.
Define $\mu P(y)= \sum_{x\in S} \mu(x) P(x,y)$ and similarly $\nu P $. Then we have
\[
%\label{relenin}
H(\mu P|\nu P) \leq H(\mu |\nu).
\]
\end{lemma}
\begin{proof}
Define $\mu_{12} (x,y)= \mu(x) P(x,y)$ and similarly $\nu_{12}(x,y)= \nu (x) P(x,y)$.
These define two joint distributions of a random variable $(X,Y)$ on $S\times S$.
Then the first marginals of $\mu_{12}, \nu_{12}$ are $\mu$, resp.\ $\nu$, and the second
marginals are $\mu P$, resp.\ $\nu P$.
Moreover, because $\sum_{y\in S} P(x,y)=1$, we get
\begin{align*}
H(\mu_{12}|\nu_{12})
&= \sum_{(x,y)\in\, S\times S} \mu(x)P(x,y) \log \frac{\mu(x)P(x,y)}{\nu(x)P(x,y)}= \sum_{x\in S} \mu(x) \log \frac{\mu(x)}{\nu(x)}\\
& = H(\mu|\nu).
\end{align*}
Therefore, by the chain rule for relative entropy (see e.g. Lemma 4.18 in \cite{rh}) we obtain
\[
H(\mu|\nu)= H(\mu P|\nu P) + D
\]
where $D$ is the conditional divergence of $X$ ``knowing'' $Y$, i.e.,
\[
D= \sum_{y\in S} \mu P(y) \sum_{x\in S} \mu_{12}(x|y)\log\frac{\mu_{12}(x|y)}{\nu_{12}(x|y)}.
\]
Because $D$ is non-negative, we obtain the desired inequality.
\end{proof}

%%% subSECTION
\subsection{Dynamics: definitions and basic inequalities}\label{subsec:dyn}

\subsubsection{Dynamics and generator}

The basic question we are interested in is how the inequality $\gcb C$ is affected by applying a Markovian dynamics to the
probability measure $\mu$. For this dynamics, we consider spin-flip dynamics with flip rates $c(i, \si)$ at site $i\in\Zd$
satisfying the following assumptions.

\subsubsection*{Condition A:}
\begin{enumerate}
\item Strict positivity: $\inf_{i\in\Zd,\,\si\in\Omega} c(i,\si)>0$.
\item Locality:
\[
%\label{ligam}
\sup_{i\in\Zd}\sum_{j\in\Zd}\,\sup_{\si\in\Omega} \big(c(i, \si^j)-c(i, \si)\big)<+\infty.
\]
\end{enumerate}
This condition ensures existence of the dynamics with generator $L$ defined below in \eqref{genflip}.

In section \ref{sec:TEGCB} we will consider weakly interacting dynamics and need more stringent conditions:

\subsection*{Condition C:}

\begin{enumerate}
\item Strict positivity: $\inf_{i\in\Zd,\,\si\in\Omega}  c(i,\si)>0$.
\item Finite-range property: There exists $R>0$ such that $c(i,\si)$ depends only on $\si_j$, for $j$ such that $|j-i|\leq R$.
\end{enumerate}
If $c(i,\si)=c(0, \tau_i\si)$, $\si\in\Omega$, $i\in\Zd$, then we say that the flip rates are translation invariant
where we remind the notation $(\tau_i\si)_j=\si_{i+j}$.

The dynamics is defined via the Markov pre-generator $L$ acting on local functions via
\be\label{genflip}
L f(\si)= \sum_{i\in \Zd} c(i,\si) \big(f(\si^i)-f(\si)\big).
\ee

As proved in \cite[Chapter 1]{ligg}, under Condition A, the closure of $L$
(in $\caC(\Omega)$ equipped with the supremum norm) generates a unique Feller process.
This process generated by $L$
is denoted $\sipr$, and $\si_i(t)$ denotes the spin at time $t$ at lattice site $i$.
We denote $\E_\si$ expectation in the process $\sipr $ starting from $\si$, and $\pee_\si$ the corresponding path-space measure.
We denote the semigroup
$S(t) f(\si)= \E_\si [f(\si(t))]$, which acts as a Markov semigroup of contractions on $\caC(\Omega)$.
Via duality, $S(t)$ acts on probability measures, and for $\mu$ a probability measure on $\Omega$, we
denote by $\mu S(t)$ the time-evolved measure, determined by the equation
\[
\int f \dd\mu S(t)= \int S(t) f \dd\mu.
\]

We also introduce the non-linear semigroup $V(t)f = \log S(t) \e^{f}$, which is a family of non-linear operators satisfying the semigroup property, i.e.,
$V(t+s)= V(t) V(s)$, $s,t\geq 0$. This non-linear semigroup appears naturally in the context of time-evolution of the Gaussian concentration bound.

Finally, notice that
\be\label{vtdel}
\left(\e^{V(t) f}\right)(\si)= \int \e^{f(\xi)} \delta_\si S(t) (\dd\xi)
\ee
whereas
\be\label{stdel}
S(t) f(\si)= \int {f(\xi)}\, \delta_\si S(t) (\dd\xi).
\ee
\subsubsection{Some basic facts for spin-flip dynamics}
In the study of existence and ergodicity properties of the Markovian dynamics $\{\si(t): t\geq 0\}$ an important role is played by the matrix
indexed by sites $i,j\in\Zd$ and defined by
\[
%\label{ligam}
\Gamma_{ij} =\sup_{\si\in\Omega} \big(c(i, \si^j)-c(i, \si)\big).
\]
We have the pointwise estimate (see \cite[Chapter 1]{ligg})
\[
%\label{}
\delta_i S(t) f\leq (\e^{t\Gamma} \delta f)_i,\,i\in\Zd, \,t\geq 0
\]
where $\e^{t\Gamma} \delta f$ denotes the bounded operator (in $\ell^1(\Zd)$) $\e^{t\Gamma}$ working on the ``column vector'' $\delta f$.
If the rates are translation invariant, i.e.,  then we have $\Gamma_{ij}= \gamma(j-i)$, i.e., $\Gamma$ acts as a convolution operator:
\[
(\Gamma \delta f)_i= (\gamma* \delta f)_i= \sum_{j\in\Zd} \gamma(i-j)\, \delta_j f
\]
and as a consequence
\[
(\e^{t\Gamma} \delta f)_i= \sum_{j\in\Zd} \gamma_t(i-j) \,\delta_j f.
\]
The so-called uniform ergodic regime, or ``$M<\varepsilon$ regime'' (see \cite{ligg}),  is the regime where the dynamics admits a unique invariant measure to which every initial measure converges
exponentially fast in the course of time. In that case there exists $\alpha>0$ such that
\be\label{mlesse}
\|\delta S(t) f\|_2^2 \leq \|\e^{t\Gamma} \delta f\|_2^2 \leq \e^{-\alpha t}\| \delta f\|_2^2
\ee
see \cite[Theorem 3.3]{ccrdyn}. In general, for a spin-flip dynamics generated by \eqref{genflip}, we have that $\Gamma$ is a bounded operator in $\ell^2(\Zd)$, i.e.,
\be\label{boundedde}
\|\delta S(t) f\|_2^2\leq K(t) \|\delta f\|_2^2
\ee
for some time-dependent constant $K(t)>0$.
Finally, we mention a useful fact about the relative entropy density.
Using the elementary Lemma \ref{relentlem}, and finite-volume approximations, one obtains the following implication for a translation invariant spin-flip dynamics
with rates satisfying condition A
\[
h(\nu|\mu)=0 \quad \Rightarrow\quad h\big(\nu S(t)\big|\,\mu S(t)\big)=0,\;\forall t>0.
\]
This will be used later on, in Section \ref{sec:nogo}.

%%% SECTION
\section{Time evolution of the Gaussian concentration bound}\label{sec:TEGCB}

In this section we show conservation of the Gaussian concentration bound under weakly interacting spin-flip dynamics, i.e., dynamics
sufficiently close to independent spin-flip dynamics.

More precisely
if we start the process $\{\si(t): t\geq 0\}$ from a probability measure $\mu$ satisfying $\gcb{C}$, then we are interested in the following questions:
\begin{enumerate}
\item Is it the case that under the time evolution $\sipr$, the time-evolved measure
$\mu S(t)$ still satisfies $\gcb{C_t}$, and if yes, how does the constant $C_t$ evolve?
\item
If the dynamics admits a unique stationary measure $\nu$, does this measure satisfy $\gcb{C}$?
\end{enumerate}

%%% subSECTION
\subsection{A general result and conservation of GCB for independent dynamics}
%\color{red}
We start  with the following general result which states that if the Gaussian concentration bound  holds at time $t>0$ when
starting from a Dirac measure $\delta_\si$ with a constant that does not depend on $\si$, then the Gaussian concentration bound
holds at time $t>0$ when started from any initial measure satisfying the Gaussian concentration bound.
%\color{black}
\begin{theorem}\label{simple}
Let $\{\si(t), t\geq 0\}$ be such that for all $\si\in\Omega$ the probability measure $\delta_\si S(t)$ satisfies $\gcb{D_t}$ where the constant
$D_t$ does not depend on $\si$. Let $\mu$ be a probability measure satisfying $\gcb{C_\mu}$.
Then, for all local functions $f$ we have
\be\label{usefulest}
\log\int \e^{f-\int f \dd\mu S(t)} \dd\mu S(t)\leq D_t\|\delta f\|_2^2  + C_\mu  \|\delta (S(t) f) \|_2^2.
\ee
As a consequence, we obtain the following results:
\begin{enumerate}
\item $\mu S(t)$ satisfies $\gcb{C(\mu, t)}$ with $C(\mu, t)\leq  D_t + K(t) C_\mu$, where $K(t)$ is defined in \eqref{boundedde}.
\item In the uniformly ergodic case \textup{(}$M<\varepsilon$ regime, cf.\ \eqref{mlesse}\textup{)},
there exists\\
$\alpha>0$ such that $\mu S(t)$ satisfies $\gcb{C(\mu, t)}$ with
\[
C(\mu, t)\leq D_t + C_\mu \e^{-\alpha t}.
\]
If furthermore, $\sup_t D_t<\infty$, then also the unique stationary measure $\nu$ satisfies $\gcb{C_\nu}$ with $C_\nu\leq \sup_t D_t<+\infty$.
\end{enumerate}
\end{theorem}
\begin{proof}
Start from the left-hand side of \eqref{usefulest}. Use that \eqref{vtdel}, \eqref{stdel}  to rewrite
\begin{align*}
& \int \e^{f-\int f \dd\mu S(t)} \dd\mu S(t) \\
& =\left(\int \left(S(t) \e^f\right) (\si)\dd\mu(\si) \right)\e^{-\int f \dd\mu S(t)} \\
&= \int \left[\left(\int \e^{f(\xi)- \int f(\zeta)\, \delta_\si S(t) (d\zeta)  } \delta_\si S(t) (\dd\xi)\right)\e^{ S(t) f(\si) - \int S(t) f (\zeta) \dd\mu(\zeta) }\right] \dd\mu(\si) \\
&\leq \e^{ D_t \|\delta f\|_2^2}\int \e^{ S(t) f(\si) - \int S(t) f (\zeta) \dd\mu(\zeta) } \dd\mu(\si) \\
& \leq \e^{ D_t \|\delta f\|_2^2} \e^{C_\mu \|\delta S(t) f\|_2^2}.
\end{align*}
In the two last steps we first used that $\delta_\si S(t)$ satisfies $\gcb{D_t}$, i.e., we have the inequality
\[
\int \e^{f(\xi)- \int f(\xi)\, \delta_\si S(t) (\dd\xi)  } \delta_\si S(t) (\dd\xi)\leq \e^{ D_t \|\delta f\|_2^2}
\]
for all $\si$.
Second, we used  the fact that $\mu$ satisfies  $\gcb {C_\mu}$. The consequences (1) and (2) now follow immediately.
\end{proof}
The following corollary shows that for independent spin-flip dynamics, Gaussian concentration is conserved.
\begin{corollary}\label{prodcol}
Assume that in the process $\sipr$ the coordinates $\{\si_i(t): t\geq 0\}$ evolve independently.
If $\mu$ satisfies $\gcb{C_\mu}$, then there exists $\alpha>0$ such that at any later time, $\mu S(t)$ satisfies $\gcb{C(\mu, t)}$, with
\be\label{coco}
C(\mu, t)= \e^{-\alpha t}C_\mu + D_t
\ee
with $\sup_t D_t <+\infty$.
\end{corollary}
\begin{proof}
First notice that if $\pee$ is a product measure on $\{-1,1\}^{\Zd}$ then $\pee$ satisfies $\gcb{C}$ with
a constant $C$ that is not depending on the marginal distributions, see \cite{blm}.
For independent spin-flip dynamics, $\delta_\si S(t)$ is a product measure. Therefore, for that case, the assumption of Theorem
\ref{simple} is satisfied, with $D_t$ uniformly bounded as a function of $t$. Furthermore, because the flip rates are assumed to be bounded from below, the process
$\sipr$ is uniformly ergodic, and as a  consequence we obtain \eqref{coco}.
\end{proof}

%%% subSECTION
\subsection{Weakly interacting spin-flip dynamics}\label{subsec:wisfd}
%\label{section4}

The result for independent spin-flip dynamics (i.e., Corollary \ref{prodcol}) can be generalized to a setting of weakly interacting dynamics, which was studied
before in \cite{mn} in the context of time-evolution of Gibbs measures. The setting is such that the rates are sufficiently close
to the rates of independent rate $1$ spin-flip dynamics, such that a space-time cluster expansion
can be set up. In particular, these conditions imply that there exists a unique invariant measure which is a Gibbs measure in the Dobrushin uniqueness regime.

More precisely, the assumptions on the rates are those of condition C, with one extra assumption forcing the rates to be close to a constant:
\be\label{extra}
c(i,\si)= 1+ \varepsilon(i, \si),\, \text{with} \, \sup_{\si\in\Omega} |\varepsilon(i, \si)| <\varepsilon_0
\ee
where $\varepsilon_0\in (0,1)$ is a constant depending on the dimension, specified in \cite{mn}.

The important implication of the space-time cluster expansion developed in \cite{mn} which we need in our context is the following.
The measure $\delta_\si S(t)$ is a Gibbs measure which is in the Dobrushin uniqueness regime, uniformly in $t>0$ and $\si$.
More precisely, $\delta_\si S(t)$ is a Gibbs measure with  uniformly absolutely summable potential $U^t_\si$ satisfying

\be\label{dobnorm}
\sup_{i\in\Zd}\sum_{\substack{A\in \caS \\ A\ni i}} |A|\sup_{\si,\eta\,\in\,\Omega,t\geq 0} \big|U^t_\si(A, \eta)\big| <1.
\ee
More precisely, in \cite{mn} an exponential norm
\[
\sup_{i\in\Zd}\sum_{\substack{A\in \caS \\ A\ni i}} \e^{a|A|} \sup_{\si,\eta\,\in\,\Omega,t\geq 0} \big|U^t_\si(A, \eta)\big|
\]
where $a>0$ is small enough, is shown to be finite, and going to zero when $\varepsilon_0\to 0$, which is stronger than \eqref{dobnorm}.

Using Theorem \ref{simple}, combined with Theorem \ref{kuelthm}, we obtain the following result.
\begin{theorem}
%\label{hight}
Let $\sipr$ be a spin-flip dynamics satisfying the conditions C, and the extra weak interaction condition \eqref{extra}.
Then we have
\begin{enumerate}
\item If $\mu$ satisfies $\gcb{C_\mu}$, then there exists $C(\mu, t)<\infty$ such that $\mu S(t)$ satisfies $\gcb{ C(\mu, t)}$.
\item The unique stationary measure $\nu$ satisfies $\gcb{C_\nu}$ for some $C_\nu<\infty$.
\end{enumerate}
\end{theorem}

%%% SECTION
\section{No-go from low-temperature Gibbs measures to Gaussian concentration bound}\label{sec:nogo}

In this section we consider a complementary regime, i.e., starting from an initial distribution where GCB is not satisfied, such as a translation-invariant Gibbs measure in the non-uniqueness regime.
We prove that it is impossible to go from such a Gibbs measure in the non-uniqueness regime towards a probability measure
which satisfies $\gcb{C}$ in finite time. One can interpret this result as the fact that one cannot acquire in finite time strong ``high-temperature'' properties from a low-temperature initial state.
%\color{red}
We prove this result first for finite-range spin-flip dynamics, and then extend to infinite range under appropriate conditions.
%\color{black}

We start with an abstract ``non-degeneracy'' condition on the Markov semigroup.
\begin{definition}[Non-degenerate Markov semigroup]
\label{nondeg}
\leavevmode\\
We say that the Markov semigroup $(S_t)_{t\geq 0}$ of a spin-flip dynamics is non-degenerate if for every pair of probability measures $\mu\not= \nu$, we have $\mu S(t)\not= \nu S(t)$ for all $t>0$.
\end{definition}
Then we have the following general result which shows that under the evolution of a non-degenerate semigroup one cannot go
from ``low temperature'' to ``high temperature'' in finite time.
\begin{theorem}
Let $\mu^+\not= \mu^-$ denote two translation-invariant Gibbs measures for the same
translation-invariant potential. Assume that the Markov semigroup is non-degenerate.
Then for all $t>0$, $\mu^+ S(t)$ cannot satisfy $\gcb{C}$.
\end{theorem}
\begin{proof}
Because $\mu^+\not= \mu^-$ are two translation-invariant Gibbs measures for the same
translation-invariant potential, we conclude that $h(\mu^- |\mu^+)=0$ and, as a consequence,
$h(\mu^- S(t) |\mu^+ S(t))=0$, for all $t>0$. By non-degeneracy, $\mu^- S(t)\not= \mu^+ S(t)$.
By \cite{moles}, we have that if $\mu^+ S(t)$ satisfies $\gcb{C}$, then
for all $\nu$ translation invariant $h_* (\nu|\mu S(t))>0$, which contradicts $h(\mu^- S(t) |\mu^+ S(t))=0$.
\end{proof}

The following lemma shows that independent spin-flip is non-degenerate.
\begin{lemma}
%\label{indfliplem}
Let $\mu, \nu$ be two different probability measures on $\Omega$. If $S(t)$ denotes the
semigroup of independent rate one spin-flip dynamics, then at any later time $t>0$, $\mu S(t)\not= \nu S(t)$.
\end{lemma}
\begin{proof}
Define, for $A\in \caS$, $\si_{\!\mathsmaller{A}}=\prod_{i\in A} \si_i$. Then we have $L \si_{\!\mathsmaller{A}}= -2|A| \si_{\!\mathsmaller{A}}$ and as a consequence,
\be\label{nodeg}
S(t)\, \si_{\!\mathsmaller{A}} = \e^{-2|A| t} \si_{\!\mathsmaller{A}}.
\ee
If $\mu S(t)= \nu S(t)$ for some $t>0$ then it follows from \eqref{nodeg} that
\[
\e^{-2|A| t} \int\si_{\!\mathsmaller{A}} \dd\mu= \e^{-2|A| t}\int \si_{\!\mathsmaller{A}} \dd\nu
\]
and therefore $\int\si_{\!\mathsmaller{A}} \dd\mu = \int \si_{\!\mathsmaller{A}} \dd\nu$.
Because linear combinations of the functions $\si_A$ are uniformly dense in $\caC(\Omega)$, we conclude that $\mu=\nu$, which leads to a contradiction.
\end{proof}

In the next subsection, we use analyticity arguments to show non-degeneracy for
general translation-invariant finite-range spin-flip dynamics.

%%% subSECTION
\subsection{Analyticity and non-degeneracy of local spin-flip dynamics}

In this section we show that for general finite-range translation-invariant spin-flip dynamics, for $\mu$ a probability measure on $\Omega$, and for a (uniformly) dense set of continuous functions $f$ the map
$t\mapsto \int S(t) f\dd\mu$ can be analytically extended to a strip in the complex plane of which the width does not depend on $\mu$.
This implies non-degeneracy in the sense of Definition \ref{nondeg}. We start with setting up the necessary notation.

We remind the notation $\si_{\!\mathsmaller{B}}=\prod_{i\in B} \si_i$ for $B$ a finite subset of $\Zd$.
For a finite set $B\subset \Zd$ we define the associated translation-invariant operator
\[
%\label{lb}
L_B= \sum_{i\in\Zd}\si_{\!\mathsmaller{B+i}} \nabla_{\!i}.
\]
In case $B=\emptyset$ we make the convention $\si_{\!\mathsmaller{B}}=1$, i.e., $L_\emptyset= \sum_{i\in\Zd} \nabla_{\!i}$ is the generator
of rate $1$ independent spin flips.

A general finite-range translation-invariant spin-flip generator can then be written in terms of these ``building block'' operators
as follows
\be\label{locca}
\loc_{\caB}:=\sum_{B\in\caB} \lambda(B) L_B
\ee
where $\caB$ is a finite collection of finite subsets of $\Zd$, and where $\lambda:\caB\to\R$. For notational simplicity, we suppressed
the dependence on the coefficient $\lambda(\cdot)$ in \eqref{locca}.
In the following lemma we produce a uniform estimate for
$L_{B_n}L_{B_{n-1}}\cdots L_{B_1} \si_{\!\mathsmaller{A}}$.
\begin{lemma}
We have the uniform estimate
\begin{align}
\nonumber
& \|L_{B_n}L_{B_{n-1}}\cdots L_{B_1} \si_{\!\mathsmaller{A}}\|_\infty\\
\label{uniest}
& \leq 2^n|A|(|A|+|B_1|)(|A|+|B_1|+|B_2|)\cdots (|A|+|B_1|+\cdots +|B_{n-1}|).
\end{align}
\end{lemma}
\begin{proof}
First notice that the bound holds when $A=\emptyset$ because  in that case $L_{B_n}L_{B_{n-1}}\cdots L_{B_1} \si_{\!\mathsmaller{A}}=0$.
So we consider $A\not=\emptyset$.
Let us first deal with $n=1$.
Notice that
\[
%\label{nablaisia}
\nabla_{\!i} \,\si_{\!\mathsmaller{A}}= -2\,\si_{\!\mathsmaller{A}} \1(i\in A)
\]
where $\1(\cdot)$ denotes the indicator function.
Next notice that $\si_{\mathsmaller{G}}\si_{\mathsmaller{F}}=\si_{\mathsmaller{G\Delta F}}$ for $G,F$ finite subsets of $\Zd$ and
$G\Delta F=(G\cap F^c)\cup (F\cap G^c)$ the symmetric difference. Then we compute
\[
L_{B_1} \si_{\!\mathsmaller{A}}= -2\,\sum_{i\,\in A} \si_{(B_1+i)\Delta A}.
\]
As a consequence
\[
\| L_{B_1} \si_{\!\mathsmaller{A}}\|_\infty \leq 2\, |A|.
\]
Let us denote for $n$ sets $C_1, \ldots, C_n$
\[
\Delta_{i=1}^n C_i= C_1\Delta C_2\Delta\cdots \Delta C_n.
\]
Then, by iteration, using $\|\si_{\!\mathsmaller{A}}\|_\infty=1$, we obtain
\begin{align*}
& L_{B_n}L_{B_{n-1}}\cdots L_{B_1} \si_{\!\mathsmaller{A}}\\
&  \quad = (-2)^n \sum_{i_1\in A} \;\sum_{i_2\in (B_1+i_1)\Delta A}\ldots\sum_{i_n\in A\Delta(\Delta_{k=1}^{n-1} (B_k+ i_k))}
\si_{A\Delta(\Delta_{k=1}^{n} (B_k+ i_k))}.
\end{align*}
Now use that $|C\Delta D|\leq |C|+ |D|$, and $\|\si_{\!A\Delta(\Delta_{k=1}^{n} (B_k+ i_k))}\|_\infty=1$, to further estimate
\begin{align*}
& \left\|\sum_{i_1\in A}\, \sum_{i_2\in (B_1+i_1)\Delta A}\ldots\sum_{i_n\in A\Delta(\Delta_{k=1}^{n-1} (B_k+ i_k))}
\si_{\!A\Delta(\Delta_{k=1}^{n} (B_k+ i_k))}\right\|_\infty \\
&\leq |A| (|A|+ |B_1|)\cdots (|A| + |B_1|+ \cdots+ |B_{n-1}|).
\end{align*}
The lemma is proved.
\end{proof}

We can then estimate $\loc_{\caB}^n\, \si_{\!\mathsmaller{A}}$.
\begin{lemma}\label{loclem}
Let $\caB\in 2^{\Zd}$ denote a finite set consisting of finite subsets of $\Zd$,
and let $\loc_\caB$ as in \eqref{locca}.
Denote $K:=\max_{B\in \caB} |B|$ and $M:=\max_{B\in \caB} |\lambda(B)|$. Then we have
\be\label{loccast}
\|(\loc_\caB^n\, \si_{\!\mathsmaller{A}})\|_\infty\leq 2^n M^n|\caB|^n(|A|+K)^n n!.
\ee
As a consequence
\[
\sum_{n=0}^\infty \frac{t^n}{n!}\left(\loc_\caB^n\, \si_{\!\mathsmaller{A}}\right)
\]
is a uniformly convergent series for $t<t_0$ with $t_0= \frac{1}{2 M|\caB|(|A|+K)}$.
\end{lemma}
\begin{proof}
We have
\[
\loc_\caB^n\, \si_{\!\mathsmaller{A}}= \sum_{B_n\in \caB}\ldots\sum_{B_1\in \caB}\left(\prod_{i=1}^n \lambda(B_i) \right)L_{B_n}\cdots L_{B_1}( \si_{\!\mathsmaller{A}}).
\]
The result then follows via \eqref{uniest} using that $|B|\leq K$ for $B\in \caB$ via the inequality
\begin{align*}
& |A| (|A|+ |B_1|)\cdots (|A| + |B_1|+ \cdots+ |B_{n-1}|)\\
& \leq  |A| (|A|+ K)\cdots (|A| + (n-1)K))\leq (|A|+K)^{n} n!.
\end{align*}
The consequence is immediate from \eqref{loccast}.
\end{proof}
\begin{proposition}\label{anal}
Let $\loc_{\caB}$ denote a finite-range translation-invariant spin-flip generator as in \eqref{locca}.
The set of analytic vectors is uniformly dense in the set of continuous functions.
\end{proposition}
\begin{proof}
The set of analytic vectors is by definition the set of functions such that there exists $t>0$ such that
\[
\sum_{n=0}^\infty \frac{t^n}{n!}\|\loc_\caB^n f\|_\infty
\]
is a convergent series. Let us denote by $\caA$ the set of analytic vectors. Notice that $\caA$ is a vector space.

By Lemma \ref{loclem} it follows that $\si_A\in\caA$ for all finite $A\subset\Zd$.
As a consequence, $\caA$
contains all local functions and as we saw before, the set of local functions is uniformly dense in $\caC(\Omega)$.
\end{proof}
\begin{proposition}\label{analprop}
Let $\mu$ and $\nu$ denote two probability measures on the configuration space $\Omega$.
Let $\loc_\caB$ denote the generator of a translation-invariant finite-range spin-flip dynamics as in \eqref{locca}.
Let $S(t)$ denote the corresponding semigroup. Let $\caA$ denote the set of analytic vectors.
Then for every $f\in \caA$, the map
\[
\psi_f (t): t\mapsto  \int S(t) f \dd\mu-\int S(t) f\dd\nu
\]
extends analytically to the set
\[
\Sigma_r:=\{ z\in \mathbb{C}: \mathrm{dist}(z,\R^+)\leq  r\}
\]
for some $r>0$ which depends on $f$ \textup{(}but not on $\mu, \nu$\textup{)}.
\end{proposition}
\begin{proof}
By assumption, there exists $r>0$ such that
\[
\sum_{n=0}^\infty \frac{t^n}{n!}\|\loc_\caB^n f\|_\infty
\]
converges
for $t\leq r$, which implies that $\psi_f (z) $ can be extended analytically in
\[
B(0,r)= \{z\in \mathbb{C}: |z|\leq r\}\subset \mathbb{C}.
\]
Now notice that the same holds when we replace $f$ by $S(s)f$, by the contraction property:
\[
\sum_{n=0}^\infty \frac{t^n}{n!}\|\loc_\caB^n(S(s)f)\|_\infty= \sum_{n=0}^\infty \frac{t^n}{n!}\|S(s)(\loc_\caB^n f)\|_\infty
\leq \sum_{n=0}^\infty \frac{t^n}{n!}\|\loc_\caB^n f\|_\infty.
\]
More precisely, for all $s$, $\psi_{S(s)f} (\cdot)$ can be extended analytically in
\[
B(0,r)= \{z\in \mathbb{C}: |z|\leq r\}\subset \mathbb{C}
\]
where $r$ does not depend on $s$.
This implies the statement of the proposition, because, via the semigroup property
\[
\psi_{S(s) f} (t)= \psi_f (s+t).
\]
The proof is finished.
\end{proof}
\begin{corollary}
Let $\mu$ and $\nu$ denote two probability measures on the configuration space $\Omega$.
Let $\loc_\caB$ denote the generator of a translation-invariant finite-range spin-flip dynamics as in \eqref{locca}.
Let $S(t)$ denote the corresponding semigroup. If $\mu\not=\nu$ then $\mu S(t)\not=\nu S(t)$  for all $t>0$.
\end{corollary}
\begin{proof}
Assume on the contrary that $\mu S(t)= \nu S(t)$ for some $t>0$, then by the
semigroup property $\mu S(s)= \nu S(s)$ for all $s\geq t$. Let $f\in \caA$ be an analytic vector such
that $\int f \dd\mu\not=\int f \dd\nu$. Then it follows that
the function $\psi_f(s)=\int S(s) f \dd\mu-\int S(s) f\dd\nu$ satisfies $\psi_f(0)\not= 0$. On the other hand, because
$\mu S(s)= \nu S(s)$ for all $s\geq t$, it follows $\psi_f(s)=0$ for all $s\geq t$. This contradicts the analyticity of
$\psi_f$.
\end{proof}

%%% subSECTION
\subsection{Generalization to a class of infinite-range dynamics}

%\textcolor{red}{Donner un exemple. Long-range Ising ?}

The assumption of finite range for the translation-invariant flip rates can be replaced by an appropriate decay condition on the
rates. This is specified below. We assume now that the generator is of the form
\[
%\label{infrange}
\caL_\caB= \sum_{B\in\caB} \lambda(B)\, \si_{\!\mathsmaller{B}} L_B
\]
where as before $L_B=\sum_{i\in\Zd} \si_{\!\mathsmaller{B+i}}\nabla_{\!i}$.
We assume now that $\caB$ is an infinite set of finite subsets of $\Zd$ and that we have the bound
\be\label{sumbo}
\sum_{B\in \caB: |B|=k} |\lambda(B)| \leq c\, \psi(k)
\ee
where $c\in (0,+\infty)$ is a constant and where $\psi(k)$ is a positive measure on the natural numbers such that for some $u>0$
\be\label{flu}
\sum_{k=0}^{\infty} \e^{u k}\psi(k)=F(u)<+\infty.
\ee

In the following lemma we obtain a bound which allows us to estimate $\|\caL_\caB^n\, \si_{\!\mathsmaller{A}}\|_\infty$.
\begin{lemma}
Let $\psi$ be a positive  measure on $\N$ such  that \eqref{flu} holds for some $u>0$.
Then for any positive integer $n$ we have
\[
\sum_{0\,\le \,k_{1},\,\ldots\,,\,0\,\le\, k_{n}}\;\prod_{j=1}^{n}\left(1+\sum_{\ell=1}^{j}k_{\ell}\right)
\;\prod_{m=1}^{n}\psi\big(k_{m}\big)\le  \e^{u} n!\;u^{-n}F(u)^{n}.
\]
\end{lemma}
\begin{proof}
We have
\begin{align*}
& \sum_{0\,\le\, k_{1},\,\ldots\,,\,0\,\le\, k_{n}}\;\prod_{j=1}^{n}\left(1+\sum_{\ell=1}^{j}k_{\ell}\right)
\;\prod_{m=1}^{n}\psi\big(k_{m}\big)\\
& \le \sum_{0\,\le\, k_{1},\,\ldots\,,\,0\,\le \,k_{n}}\;\left(1+\sum_{\ell=1}^{n}k_{\ell}\right)^{n}
\;\prod_{m=1}^{n}\psi\big(k_{m}\big)
\\
&=\e^{u}\;\sum_{0\,\le\, k_{1},\,\ldots\,,\,0\,\le\, k_{n}}\;\e^{-u\,\big(1+\sum_{j=1}^{n}k_{j}\big)}\;
\left(1+\sum_{\ell=1}^{n}k_{\ell}\right)^{n}
\;\prod_{m=1}^{n}\left(\e^{u\,k_{m}}\psi\big(k_{m}\big)\right)\\
& \le \e^{u}n! \,u^{-n}F(u)^{n}
\end{align*}
where we used that $v^{n} \e^{-v}\!/n!< 1$, for all $v>0$ and $n$.
\end{proof}

We can then show that the bound of Lemma \eqref{uniest} still holds.
\begin{proposition}
Under \eqref{sumbo} and \eqref{flu}, we have the bound
\[
%\label{infbound}
\|\loc_\caB^n\, \si_{\!\mathsmaller{A}}\|_\infty\leq n!\,\kappa^n, \; n\geq 1,
\]
for some $\kappa>0$. As a consequence, local functions are analytic vectors, and the Markovian dynamics generated by $\loc_\caB$ is non-degenerate.
\end{proposition}
\begin{proof}
We estimate as in the proof of Lemma \ref{loclem}, using \eqref{flu}. Let $u>0$ be as in \eqref{flu}, and $A\not=\emptyset$. Then
\begin{align*}
& \|\loc_\caB^n\, \si_{\!\mathsmaller{A}}\|_\infty \\
& \leq
2^n\sum_{B_1\in \,\caB}\cdots \sum_{B_n\in \,\caB_n} \\
& \quad \times \;  \sum_{i_1\in A} \,\sum_{i_2\in (B_1+i_1)\Delta A}\cdots\sum_{i_n\in A\Delta(\Delta_{k=1}^{n-1} (B_k+ i_k))}
\!\|\si_{\!A\Delta(\Delta_{k=1}^{n} (B_k+ i_k))}\|_\infty \prod_{i=1}^n|\lambda(B_i)| \\
& \leq
2^n\sum_{k_1=0}^\infty\cdots \sum_{k_n=0}^\infty \sum_{B_1\in \caB, |B_1|=k_1}\cdots \sum_{B_n\in \caB, |B_n|=k_n}
|A| (|A|+k_1)\cdots \\
& \qquad \times \; (|A|+k_1+k_2+\cdots +k_n)\prod_{i=1}^n|\lambda(B_i)| \\
& \leq 2^n c^n\sum_{k_1=0}^\infty\cdots \sum_{k_n=0}^\infty
|A| (|A|+k_1)\cdots (|A|+k_1+k_2+\cdots +k_n)\prod_{i=1}^n \psi(k_i) \\
& \leq |A|^n \,2^n c^n \e^{u} n!\,u^{-n} F(u)^{n} \\
& \leq n! \,\kappa^n
\end{align*}
for some $0<\kappa<+\infty$. With this bound, we can proceed as in the proof of the finite-range case (Lemma \ref{loclem}, Propositions \ref{anal},
\ref{analprop}).
\end{proof}

%%% SECTION
\section{Uniform variance bound}\label{sec:UVB}

In this section we consider the time-dependent behavior of a weaker concentration inequality, which we call the ``uniform variance bound''. In the context of Gibbs measures, contrarily to GCB,  this inequality can still
hold in the non-uniqueness regime (for the ergodic equilibrium states), see \cite{cckr} for a proof of this inequality for the
low-temperature pure phases of the Ising model.

\begin{definition}[Uniform Variance Bound]
\leavevmode\\
We say that $\mu$ satisfies the uniform variance bound with constant $C$ (abbreviation $\uvi{C}$) if for all $f:\Omega\to\R$ continuous
\be
\label{uvic}
\E_\mu \big[(f-\E_\mu(f))^2\big]\leq C\|\delta f\|_2^2.
\ee

\end{definition}
\begin{remark}
%\color{red}
The bound \eqref{uvic} implies that
ergodic averages of the form
$(\sum_{x\in\la} \tau_x f)/|\la|$
(with $f$ a local function) have a variance which is bounded
by $C_f/{|\la|}$. In particular this excludes convex combinations of pure phases (non-ergodic states), and critical behavior (states
at the critical point). The result which we show below (Theorem \ref{fifombo}) thus shows that when started from
an initial measure satisfying \eqref{uvic},  non-ergodic or critical behavior cannot be obtained in finite time.
%\color{black}
\end{remark}

Notice that, in contrast with the Gaussian concentration bound, the inequality \eqref{uvic} is homogeneous, i.e.,
if \eqref{uvic} holds for $f$ then for all $\lambda\in \R$, it also holds for $\lambda f$. Furthermore, if
\eqref{uvic} holds for a  subset of continuous functions which is uniformly dense in $\caC(\Omega)$ (such as the set of local functions), then it
holds for all $f\in \caC(\Omega)$ by standard approximation arguments.
This implies that if we can show the validity of \eqref{uvic} for a set of functions $\caD\subset \caC(\Omega)$ such that
$\cup_{\lambda\in [0,+\infty)} \lambda \caD$ contains all local functions, we obtain the validity of \eqref{uvic} for all $f\in \caC(\Omega)$.

The following proposition shows that a weak form of Gaussian concentration is equivalent with the uniform variance bound.
\begin{definition}[Weak Gaussian Concentration Bound]
%\label{def-wgcb}
\leavevmode\\
We say that a probability measure $\mu$ satisfies the weak Gaussian concentration bound with constant $C$ if for
every $f:\Omega\to\R$ continuous there exists $\lambda_0=\lambda_0(f)>0$ such that for all $\lambda\leq \lambda_0$
\be\label{wgcb}
\E_\mu\left( \e^{\lambda(f-\E_\mu (f))}\right)\leq \e^{C\lambda^2\|\delta f\|_2^2}.
\ee
\end{definition}
\begin{proposition}\label{gcbuvilem}
A probability measure $\mu$ satisfies the weak Gaussian concentration bound with constant $C$ if and only if it satisfies the uniform variance bound.
\end{proposition}
\begin{proof}
Assume that $\mu$ satisfies the weak Gaussian concentration bound with constant $C$.
From  \eqref{wgcb} we derive, for $f:\Omega\to\R$ continuous,
\[
\var_\mu(f)= \lim_{\lambda\to 0}\frac{\E_\mu\left(\e^{\lambda(f-\E_\mu (f))}\right)-1}{\lambda^2}\leq \lim_{\lambda\to 0}\frac{\e^{C\lambda^2\|\delta f\|_2^2}-1}{\lambda^2}= C\|\delta f\|_2^2.
\]
which is the uniform variance bound.
Conversely, assume that the uniform variance bound holds, and let $f:\Omega\to\R$ be a continuous function. Then use
the elementary inequality $\e^{\lambda x} -1-\lambda x\leq \frac{\lambda^2 \mathrm{e}\,x^2}{2}$, valid for
for $0\leq\lambda x\leq 1$, together with $\e^x\geq 1+x$, to conclude that for $\lambda \leq \frac{1}{2\|f\|_\infty +1}$, we have
\[
\E_\mu\left( \e^{\lambda(f-\E_\mu (f))}\right)\leq 1+\frac{\lambda^2 \mathrm{e}\var_\mu (f)}{2}
\leq 1+ \frac{\lambda^2\e}{2}\, C\|\delta f\|_2^2 \leq \e^{\tfrac{\lambda^2\mathrm{e}}{2}C\|\delta f\|_2^2}.
\]
\end{proof}
%\color{red} Je pense qu'on a la m\^eme constante $C$ car il existe $u_0>0$ tel que $\e^u-1-u\leq u^2$, $\forall u\leq u_0$.\color{black}
The following theorem is the analogue of Theorem \ref{simple} for the uniform variance bound.
\begin{theorem}\label{simplevar}
Assume that for all $\si$, the probability measure $\delta_\si S(t)$ satisfies $\uvi{C(\si, t)}$. If $\mu$ satisfies $\uvi{C}$ and is such that $\int C(\si, t) \dd\mu(\si)<+\infty$, then also $\mu S(t)$ satisfies $\uvi{C(\mu, t)}$
with
\[
%\label{unic}
C(\mu, t)\leq  CK(t) + \int C(\si, t) \dd\mu(\si)
\]
where $K(t)$ is as in \eqref{boundedde}.
\end{theorem}
\begin{proof}
Let $f:\Omega\to\R$ be a continuous function. Then we compute, using \eqref{boundedde}:
\begin{align*}
\var_{\mu S(t)} (f)
& = \int f^2 \dd\mu S(t) - \left(\int f \dd\mu S(t)\right)^2 \\
& =\int \left(S(t) (f^2) - (S(t) f)^2\right) \dd\mu + \var_\mu (S(t) f) \\
& = \int \var_{\delta_\si S(t)} (f) \dd\mu (\si) + C \|\delta S(t) f\|_2^2 \\
& \leq \left(\int C(\si, t) \dd\mu (\si)\right) \|\delta f\|_2^2 + C K(t)^2 \|\delta f\|^2_2.
\end{align*}
The theorem is proved.
\end{proof}
\begin{corollary}\label{highuvb}
Assume that the spin-flip rates satisfy the weak interaction condition of Section \ref{subsec:wisfd}, then the dynamics conserves the uniform variance bound.
\end{corollary}
\begin{proof}
Under the weak interaction condition, $\delta_\si S(t)$ satisfies $\gcb C$ with a constant that does not depend on $\si$.
By Proposition \ref{gcbuvilem} $\delta_\si S(t)$ satisfies $\uvi{C}$ with a constant that does not depend on $\si$.
The conclusion follows from Theorem \ref{simplevar}.
\end{proof}

The following theorem shows that the high-temperature condition of corollary \ref{highuvb} is not necessary, and in fact,  the uniform variance inequality is robust under any local spin-flip dynamics, i.e.,
under the condition C of Section \ref{subsec:dyn}.
\begin{theorem}\label{fifombo}
Assume that $\mu$ satisfies the uniform variance inequality \eqref{uvic}. Let $S(t)$ denote the semigroup of a spin-flip dynamics
condition A of Section \ref{subsec:dyn}.
Then $\mu S(t) $ satisfies the uniform variance inequality for all $t>0$.
\end{theorem}
\begin{proof}
Let us denote the time-dependent quadratic form
\[
%\label{quadr}
\psi(t;f,g)= S(t)(fg)- (S(t)f) (S(t)g)
\]
as well as the usual carr\'{e} du champ quadratic form
\[
%\label{gamma}
\Gamma(f,g)= L(fg)- gLf-fLg.
\]
Notice that
\be\label{corambo}
\var_{\delta_\si S(t)} (f)= \psi(t; f,f)(\si) .
\ee
An simple explicit computation shows that
\[
%\label{spingamma}
\Gamma(f,f)= \sum_{i\in\Zd} c(i,\si) (f(\si^i)-f(\si))^2
\]
which by the boundedness of the rates implies the estimate
\be\label{gamest}
\|\Gamma(f,f)\|_\infty\leq \hat{c}\, \|\delta f\|_2^2
\ee
with $\hat{c}=\sup_{\si\in \Omega, i\in\Zd} c(i,\si)$.
%\textcolor{red}{RENVOYER aux notes de van Handel p. 41 (\`a v\'erifier).}\\
We then compute
\begin{align*}
\frac{\dd}{\dd t} (\psi(t;f,f))
&= L(S(t) f^2)- 2 S(t)f L S(t) f \\
&= L\left[S(t)f^2 - (S(t) f)(S(t)f)\right]+ 2\Gamma (S(t) f, S(t) f).
\end{align*}
As a consequence, using $\psi(0;f,f)=0$, by the variation of constants method we obtain
\[
\psi(t;f,f)=2\int_0^t S(t-s)\,\Gamma(S(s)f, S(s) f) \dd s.
\]
Therefore, using \eqref{gamest} combined with the contraction property of the semigroup, we obtain, via \eqref{boundedde}
\[
\|\psi(t;f,f)\|_\infty \leq 2\,\hat{c}\int_0^t \|\delta S(s) f\|_2^2 \dd s\leq 2\,\hat{c} \left(\int_0^t K(s)^2 \dd s\right) \|\delta f\|_2^2.
\]
Now use \eqref{corambo}
to conclude
\[
\var_{\delta_\si S(t)} (f)\leq C \|\delta f\|_2^2
\]
with $C= 2\, \hat{c} \left(\int_0^t K(s)^2 \dd s\right)$ not depending on $\si$.

Via Theorem \ref{simplevar},  we obtain the statement of the theorem.
\end{proof}
\begin{remark}
Remark that we did not use the finite range character of the spin-flip rates, neither the translation invariance.
I.e., as soon as the flip rates are uniformly bounded, and are such that the Markovian dynamics with these rates can be defined, we obtain that the uniform variance bound is conserved in the course of time.
\end{remark}

Finally, we show the analogue of Theorem \ref{simple} for more general inequalities including moment inequalities.
\begin{definition}
Let $\euH:\R\to \R$ be a convex function,  $\euJ: [0, \infty)\to\R$ a continuous increasing function, and $C>0$ a constant. Then we say that
$\mu$ satisfies the $(\euH,\euJ, C)$ inequality if for all continuous $f:\Omega\to\R$ with $\|\delta f\|_2<+\infty$ we have
\[
%\label{HJC}
\int \euH(f-\E_\mu(f)) \dd\mu\leq \euJ(C\|\delta f\|_2).
\]
\end{definition}
To fit the  examples we saw so far: we have $\uvi{a}$ corresponds to $\euH(x)=x^2, \euJ(x)=x^2, C=\sqrt{a}$, whereas $\gcb{a}$ corresponds to $\euH(x)= \e^{x}, \euJ(x)= \e^{x^2}, C=\sqrt{a}$.
More general moment inequalities correspond to $\euH(x)= |x|^p, \euJ(x)= |x|^p$.

The following theorem is then the analogue of Theorem \ref{simple} for the $(\euH,\euJ,C)$ inequality.

\begin{theorem}
Assume that $\delta_\si S(t)$ satisfies the $(\euH,\euJ,C)$ inequality with constant $C$ that does not depend on
$\si$. Then if $\mu$ satisfies the $(\euH,\euJ,C_\mu)$ inequality, so does $\mu S(t)$ for all $t>0$.
\end{theorem}

\begin{proof}
We write, using $p_t(\si, \dd\eta)$ for the transition probability measure starting from $\si$, and abbreviating $\int f \dd\mu S(t) =: \mu(t, f)$
\begin{align*}
& \int \euH(f-\mu(t,f)) \dd\mu S(t)\\
&=
\int p_t (\si, \dd\eta)\,  \euH\left( f(\eta)-\int f(\xi)\, p_t(\si, \dd\xi) + \int f(\xi)\, p_t(\si, \dd\xi)- \mu(t, f)\right)\dd\mu(\si)\\
&\leq
\frac12 \int p_t(\si, \dd\eta)\, \euH\left(2\left(f(\eta)-\int f(\xi)\, p_t(\si, \dd\xi)\right)\right)\dd\mu(\si)\\
& \qquad + \frac12 \int \euH\left(2 \left(S(t)f- \int S(t) f \dd\mu\right)\right)\dd\mu(\si) \\
&\leq
\frac12 \euJ( 2C \|\delta f\|_2) + \frac12 \euJ(2C_\mu \|\delta S(t) f\|_2)\\
&\leq \frac12 \euJ( 2C \|\delta f\|_2) + \frac12 \euJ\big(2C_\mu \sqrt{K(t)}\, \|\delta f\|_2\big)\\
&\leq
\euJ\big( \big(2C+2C_\mu \sqrt{K(t)}\, \big)\|  \delta f\|_2\big).
\end{align*}
Here in the last two steps we used \eqref{boundedde}, combined with the fact that $\euJ$ is increasing.
\end{proof}

%%% BIBLIOGRAPHY

\end{document}